\documentclass[12pt,twoside]{article}    
\usepackage{amssymb}
\usepackage{amsmath}
\usepackage{amsfonts}
\usepackage{graphicx}
\usepackage{url}

\providecommand{\U}[1]{\protect\rule{.1in}{.1in}}

\newtheorem{theorem}{Theorem}[section]
\newtheorem{remark}[theorem]{Remark}
\newtheorem{lemma}[theorem]{Lemma}

\newtheorem{example}[theorem]{Example}
\numberwithin{equation}{section}
\newenvironment{proof}[1][Proof]{\textbf{#1.} }{\ \rule{0.5em}{0.5em}}

\makeatletter \@addtoreset{equation}{section} \makeatother

\begin{document}


\pagestyle{myheadings}

\markboth{\hfill {\small AbdulRahman Al-Hussein} \hfill}{\hfill
{\small FBDSDEs with Poisson jumps} \hfill }


\thispagestyle{plain}


\begin{center}
{\large \textbf{Forward-backward doubly stochastic differential equations with Poisson jumps in infinite dimensions}} \\
\vspace{0.7cm} AbdulRahman Al-Hussein
\\
\vspace{0.2cm} {\footnotesize
{\it Department of Mathematics, College of Science, Qassim University, \\
 P.O.Box 6644, Buraydah 51452, Saudi Arabia \\ {\emph E-mail:} alhusseinqu@hotmail.com, hsien@qu.edu.sa}}
\end{center}

\vspace{0.15cm}

\begin{abstract}
In this paper, we study the existence and uniqueness of solution to a system of nonlinear fully coupled forward-backward doubly stochastic differential equations with Poisson jumps. Our work is established in infinite dimensional separable Hilbert spaces and is based on the method of time continuation.
\end{abstract}

{\bf MSC 2010:} 60H10, 60G55. \\

{\bf Keywords:} Continuation method, forward-backward doubly stochastic differential equations, monotonicity condition, Poisson process, Wiener process.

\bigskip

\section{Introduction}
Fully coupled forward-backward doubly stochastic differential equations (FBDSDEs, for short) were studied by Peng and Shi in \cite{PS} to generalize stochastic Hamiltonian systems. The equations considered in their work operate on the same Euclidean space. Later, the work was extended in \cite{ZSG} to encompass FBDSDEs in different dimensions. 

On the other hand, backward stochastic differential equations (BSDEs, shortly) in infinite dimensions were studied in \cite{Alh-2004}, \cite{Hu-Peng-90}, and \cite{Tessitore}. In particular, Al-Hussein, \cite{Alh-2004}, focused on the existence and uniqueness of the solutions of BSDEs driven by genuine $\mathcal{Q}$-Wiener processes or cylindrical Wiener processes on separable Hilbert spaces. Furthermore, Al-Hussein and Gherbal, \cite{Alh-Ghe-Rose21}, established the theory of FBDSDEs with Poisson jumps in Euclidean spaces. Recent developments in this field can be found in \cite{ZST} and the relevant references therein. In fact, in \cite{Alh-Ghe-Rose21}, three different cases were distinguished based on the dimensions of the Euclidean spaces where the solution processes evolve. However, when working more generally over infinite dimensions, these cases necessarily change. To the best of our knowledge, apart from Theorem~2 in our previous work \cite{AG}, which lacked a proof, no other work in the literature addresses such systems of FBDSDEs in infinite dimensions.

The objective of this paper is to address this research gap by providing a more comprehensive treatment of FBDSDEs in infinite dimensions. We aim to establish the theory of existence and uniqueness of the solution to fully coupled FBDSDEs with jumps (FBDSDEJ, for short) on separable real Hilbert spaces. The system we consider is given by:
\begin{eqnarray}\label{intr-syst}
\left\{
\begin{array}{ll}%
dy_{t}=b\left(  t,y_{t},Y_{t},z_{t},Z_{t},k_{t}\right)
dt+\sigma\left(  t,y_{t},Y_{t},z_{t},Z_{t},k_{t}\right)  dW_{t}%
-z_{t}\, \overleftarrow{dB}_{t}\\ \hspace{2.5in}
+\int_{\Theta}\varphi\left(  t,y_{t},Y_{t},z_{t},Z_{t},k_{t},\rho\right)
\widetilde{N}(d\rho,dt), \\ \\
dY_{t}=f\left(  t,y_{t},Y_{t},z_{t},Z_{t},k_{t}\right)  dt+g\left(  t,y_{t},Y_{t},z_{t},Z_{t},k_{t}\right)
 \overleftarrow{dB}_{t}
\\ \hspace{2.5in}
+ \, Z_{t}\, dW_{t} + \int_{\Theta}k_{t}\left(  \rho\right)  \widetilde{N}(d\rho,dt) ,
\\ \\
y_{0}=x  ,Y_{T}=h\left(  y_{T}\right) .
\end{array}
\right.
\end{eqnarray}
Here, $(W_{t})_{t\geq0}$ and
$(B_{t})_{t\geq0}$ are independent cylindrical Wiener processes on separable real Hilbert spaces $E_1$ and $E_2$ respectively, and $\tilde{N}(d\rho,dt)$ represents a compensated Poisson random measure.
The mappings $b,\sigma,\varphi,f,g$, and $h$ are defined and discussed in Section~\ref{sec3}.
The stochastic integral appearing with a backward arrow, such as $\int_0^t z_{s}\, \overleftarrow{dB}_{s}$, is a backward It\^{o} integral.

We emphasize that the system (\ref{intr-syst}) is Markovian and the imposed conditions are natural and accessible for verification,  enabling the development of concrete examples and counterexamples, as demonstrated in Examples~(\ref{ex:ex1}, \ref{ex:ex2}) in Section~\ref{sec3}. Moreover, this work is crucial for addressing stochastic control problems associated with FBDSDEJ in infinite dimensions using the maximum principle approach. Further insights on this topic can be found in our previous works, namely
 \cite{AG} and \cite{Al-G-relaxed}. Recently, our research in \cite{Alh-Ghe-022} highlighted the significant role of this paper in deriving strong solutions for FBDSDEs of McKean-Vlasov type. FBDSDEs have also shown relevance in optimal filtering problems, as observed in \cite{Bao}.
  Finally, it is worth mentioning that numerical investigations in 
  \cite{Hu-Nua-Song011} and \cite{Hu-Nua-Song017}, employing Malliavin calculus for FBSDEs and BDSDEs, respectively, have led to promising research on  a linear version of FBDSDEJ~(\ref{intr-syst}) in both finite and infinite dimensions.
\subsection{Notation and Spaces of Solutions}
 Let $(\Omega, \mathcal{F},\mathbb{P})$ be a complete probability space. Let $(W_t)_{t\in [0,T]}$ and $(B_t)_{t\in [0,T]}$ be two cylindrical Wiener process on
 separable real Hilbert spaces $E_1$ and $E_2$, respectively, and let $\eta$ be a Poisson point process taking
its values in a measurable space $( \Theta,\mathcal{B}( \Theta) ),$ where $\Theta$ is a standard Borel space.
Denote by $\Pi( d\rho) $ to the characteristic measure of $\eta$ which is assumed to be a
$\sigma$-finite measure on $( \Theta,\mathcal{B}( \Theta) ),$ by $N (d\rho,dt)
$ to the Poisson counting measure induced by $\eta$ with compensator $\Pi( d\rho)
dt$, and by
\(
\tilde{N}( d\rho,dt) = N ( d\rho,dt) -\Pi( d\rho) dt
\)
to the compensation of the jump measure $N(\cdot, \cdot)$ of $\eta .$ Then $\Pi (O) = \mathbb{E} [N(O,1)]$ for ${O\in \mathcal{B}(\Theta).}$
The three processes $W , B$ and $\eta$ are assumed to be mutually independent. Let
\begin{eqnarray*}
 && \mathcal{F}_t^W := \sigma \left(\{\ell(W_r) | \; 0 \leq r\leq t, \, \ell \in E_1^*\} \right) \vee \mathcal{N}, \\
 && \mathcal{F}_{t,T}^B := \sigma \left(\{\ell(B_r -B_T) | \; t\leq r\leq T, \, \ell \in E_2^* \}\right) \vee \mathcal{N}, \\
 && \mathcal{F}_{t}^N:=\sigma\{N_{r} | \; 0\leq r\leq t \}\vee \mathcal{N},
\end{eqnarray*}
for all $0\leq t\leq T$, where $\mathcal{N}$ consists of all the $\mathbb{P}-$null sets in $\Omega$.

It is obvious that $\{\mathcal{F}_{t}^{W}\mid t\geq 0\}$ and $\{\mathcal{F}_{t}^N\mid t\geq 0\}$ are filtrations in the sense that $\mathcal{G}_{t}\subseteq \mathcal{G}_{s}$, if $t\leq s,$ for $\mathcal{G}=\mathcal{F}^{W}, \mathcal{F}^N$, while $\{\mathcal{F}_{t,T}^{B}\mid 0\leq t\leq T\}$ is a backward filtration in the sense that $\mathcal{F}_{t,T}^{B}\supseteq \mathcal{F}_{s,T}^{B}$ if $s\geq t.$
More details about backward filtrations and backward martingales can be found in \cite{Alh-Ghe-Rose21}.

Denote by $L_2(E_2;H)$ to the space of Hilbert-Schmidt operators from $E_2$ to $H$, which is also a Hilbert space with respect to its usual norm $| \cdot |_{L_2(E_2;H)}$. Let $h$ be an $\{\mathcal{F}^{B}_{t,T}\mid 0\leq t\leq T\}$-adapted $L_2(E_2; H)$-valued process satisfying $\mathbb{E} [\int_{0}^{T}|h_{s}|^{2}ds]<\infty$. The \emph{backward It\^{o} integral} of mean square continuous such a stochastic process $h$ with respect to $B$ is given by
\begin{equation*}
\int_{\alpha}^{\beta}h_{s}\, \overleftarrow{dB}_{s}=\lim_{|\pi|\rightarrow 0}\,\overset{n}{\underset{i=1}{\mathbb{\sum}}}h(t_{i+1})(B_{t_{i+1}}-B_{t_{i}}) \quad \quad (\text{in} \; L^{2}(\Omega,\mathcal{F},\mathbb{P}; H)),
\end{equation*}
where $\pi:=\{t_{1},t_{2},\cdots,t_{n+1}\}$ is a partition of $[\alpha,\beta]$ with $\underset{1\leq i\leq n}{\max}(t_{i+1}-t_{i})\rightarrow 0.$

Next, let us denote
 $\breve{B}_{s}:=B_{T-s}-B_{T},\; 0\leq s\leq T.$ Then, obviously, $\breve{B}$ is a cylindrical Wiener process as well, and for all $0\leq s\leq T$, we have
$\mathcal{F}_{T-t,T}^{B}=\mathcal{F}_{t}^{\breve{B}}, $ and similarly for $\breve{W}_{s}:=W_{T-s}-W_{T}.$
We have also $\mathcal{F}_{T-t,T}^{\breve{N}}=\mathcal{F}_{t}^N$ for all $0\leq t\leq T.$
Therefore, if $h_{s}$ is $\mathcal{F}_{s,T}^{B}$-measurable for each $0\leq s\leq T,$ then $\breve{h}_{s}:=h_{T-s}$ is $\mathcal{F}_{s}^{\breve{B}}$-measurable for each $0\leq s\leq T,$ and if $k_{s}$ is $\mathcal{F}_{s}^{W}$-measurable for each $0\leq s\leq T,$ then $\breve{k}_{s}:=k_{T-s}$ is $\mathcal{F}_{T-s,T}^{\breve{W}}$-measurable for each $0\leq s\leq T.$

If $\{h_{t}\mid 0\leq t\leq T\}$ is an $L_2(E_2;H)$-valued stochastic process such that $h_t$ if $\mathcal{F}_{t}\triangleq\mathcal{F}_{t}^{W}\vee\mathcal{F}_{t,T}^{B}\vee\mathcal{F}_{t}^N$-measurable for each $0\leq t\leq T$ and $\mathbb{E}\, [\int_{0}^{T}|h_{t}|^{2}dt]<\infty$,  then for $\breve{h}_{s}\, (:=h_{T-s})$,  we get 

\begin{equation}\label{eq:2.3}
\int_{0}^{T-t}\breve{h}_{s}\, d\breve{B}_{s}=-\int_{t}^{T}h_{s}\, \overleftarrow{dB}_{s}, \quad 0\leq t\leq T,
\end{equation}
see \cite{Alh-Ghe-Rose21}. Hence,
\begin{equation}\label{eq:2.6}
\int_{0}^{u}h_{s}\, \overleftarrow{dB}_{s}=-\int_{T-u}^{T}\breve{h}_{s}\, d\breve{B}_{s}, \quad 0\leq u\leq T.
\end{equation}

For a separable Hilbert space $E,$ let $\mathcal{M}^{2}(0,T; E) $ denote the set of jointly measurable,
processes $\left\{   \mathcal{X}_{t},t\in\big[ 0,T  \big] \right\}
$ with values in $E$ such that $ \mathcal{X}_t$ is $\mathcal{F}_{t}$-measurable for a.e. $t \in [0,T]$ and
\[
|\mathcal{X}|^2:=\mathbb{E} \left[  {\int_{0}^{T}} \left\vert \mathcal{X}_{t}\right\vert_E ^{2}dt \right]
<\infty.
\]
Denote by $L_{\Pi}^{2}( E)$ to the set of $\mathcal{B}( \Theta)$-measurable mapping $k$ with values in $E$
such that
\[
||| k ||| := \left[ \int_{\Theta}\left\vert k(
\rho) \right\vert_{E}^{2} \, \Pi( d\rho) \right]^{\frac{1}{2}} < \infty .
\]
Let also $\mathcal{N}_{\eta}^{2}( 0,T ;E) $ be the space of all $L_{\Pi}^{2}( E)$-valued processes $\{ \mathcal{K}_t , \; t \in \big[0,T \big] \}$ that satisfy: $\mathcal{K}_t$ is $\mathcal{F}_{t}$-measurable for a.e. $t \in [0,T]$ and
\[
\mathbb{E} \left[ {\int_{0}^{T}} \int_{\Theta}\left\vert \mathcal{K}_{t}( \rho)
\right\vert_{E}^{2}\Pi( d\rho) dt \right] <\infty.
\]

Using the preceding notation we set
\begin{align*}
& \mathfrak{M}^{2} := \mathcal{M}^2(0,T;H)\times \mathcal{M}^2(0,T;H) \times \mathcal{M}^2(0,T;L_2(E_2;H)) \\
& \hspace{2.5in} \times \mathcal{M}^2(0,T;L_2(E_1;H)) \times\mathcal{N}_{\eta}^{2}\left(0,T ;H\right) .
\end{align*}
This space $\mathfrak{M}^{2}$ is evidently a Hilbert space
 with the norm $\left\Vert \cdot \right\Vert _{\mathfrak{M}^{2}}$ given by:
\begin{equation*}
\left\Vert \zeta_{\cdot}\right\Vert _{\mathfrak{M}^{2}}^{2} := \mathbb{E}\left[ \int_{0}^{T} \Big{(} \left\vert
y_{t}\right\vert ^{2} + \left\vert Y_{t}
\right\vert ^{2} + \left\Vert z_{t}\right\Vert^{2} + \left\Vert Z_{t}\right\Vert^{2} + |||
k_{t}||| ^{2} \Big{)} dt\right]  ,
\end{equation*}
for $\zeta_{\cdot}=\left(  y_{\cdot},Y_{\cdot},z_{\cdot},Z_{\cdot},k_{\cdot}  \right) \in \mathfrak{M}^{2} .$

\section{Main Results}\label{sec3}
Let us denote 
$\mathbb{H}^2:= H \times H \times L_2(E_2;H)\times L_2(E_1;H)\times L_{\Pi }^{2}( H)$. 
Consider
\begin{eqnarray*}
&& b :\Omega \times [0,T]  \times  \mathbb{H}^2  \rightarrow H ,\\
&& \sigma :\Omega \times [  0,T]  \times  \mathbb{H}^2 \rightarrow L_2(E_1;H),\\
 && \varphi :\Omega \times [  0,T]  \times \mathbb{H}^2  \times \Theta \rightarrow H ,\\
 && f  :\Omega \times [  0,T]  \times \mathbb{H}^2 \rightarrow
 H ,\\
&& g :\Omega \times [  0,T]  \times \mathbb{H}^2\rightarrow L_2(E_2;H),\\
 && h :\Omega \times H \rightarrow H ,
\end{eqnarray*}
as mappings with properties to be given shortly.
Consider also the FBDSDEJ:
\begin{eqnarray}\label{eq:3.1}
\left\{
\begin{array}{ll}%
dy_{t}=b\left(  t,y_{t},Y_{t},z_{t},Z_{t},k_{t}\right)
dt+\sigma\left(  t,y_{t},Y_{t},z_{t},Z_{t},k_{t}\right)  dW_{t}%
-z_{t}\, \overleftarrow{dB}_{t}\\ \hspace{2.5in}
+\int_{\Theta}\varphi\left(  t,y_{t},Y_{t},z_{t},Z_{t},k_{t},\rho\right)
\widetilde{N}(d\rho,dt), \\ \\
dY_{t}=f\left(  t,y_{t},Y_{t},z_{t},Z_{t},k_{t}\right)  dt+g\left(  t,y_{t},Y_{t},z_{t},Z_{t},k_{t}\right)
\overleftarrow{dB}_{t}
\\ \hspace{2.5in}
+\, Z_{t}\, dW_{t}
+\int_{\Theta}k_{t}\left(  \rho\right)  \widetilde{N}(d\rho,dt) ,
\\ \\
y_{0}=x  ,Y_{T}=h\left(  y_{T}\right) .
\end{array}
\right.
\end{eqnarray}
Here $x$ is a fixed element of $H$.

\smallskip

A \emph{solution} of (\ref{eq:3.1}) consists of stochastic processes $(y,Y,z,Z,k)$, such that for each $t\in[0,T]$, we have $a.s.:$
\begin{eqnarray}\label{eq:3.2}
\left\{
\begin{array}{ll}%
y_{t}=x+\int_{0}^{t}b\left(  s,y_{s},Y_{s},z_{s},Z_{s},k_{s}\right)
ds+\int_{0}^{t}\sigma\left(  s,y_{s},Y_{s},z_{s},Z_{s},k_{s}\right)  dW_{s}
\\ \hspace{1in}
-\, \int_{0}^{t}z_{s}\, \overleftarrow{dB}_{s} +\int_{0}^{t}\int_{\Theta}\varphi\left(  s,y_{s},Y_{s},z_{s},Z_{s},k_{s},\rho\right)
\widetilde{N}(d\rho,ds),\\ \\
Y_{t}=h\left(  y_{T}\right)-\int_{t}^{T}f\left(  s,y_{s},Y_{s},z_{s},Z_{s},k_{s}\right)  ds-\int_{t}^{T}Z_{s}\, dW_{s}\\ \hspace{1in}
-\int_{t}^{T}g\left(  s,y_{s},Y_{s},z_{s},Z_{s},k_{s}\right)  \overleftarrow{dB}_{s}
-\int_{t}^{T}\int_{\Theta}k_{s}\left(  \rho\right)  \widetilde{N}(d\rho,ds).
\end{array}
\right.
\end{eqnarray}

To simplify the notation in what follows, we let
$$\upsilon  =\left(  y,Y,z,Z,k\right)  , \; A\left(  t,\upsilon\right)    :=\left(    f,  b,  g,  \sigma,  \varphi\right)  \left(
t,\upsilon\right)  ,$$
$$\Arrowvert \upsilon \Arrowvert^2 := \left\vert
y \right\vert ^{2} + \left\vert Y 
\right\vert ^{2} + \left\Vert z \right\Vert^{2} + \left\Vert Z \right\Vert^{2} + |||
k ||| ^{2} , $$
$$\left\langle A,\upsilon\right\rangle  :=\left\langle
y, f\right\rangle_H +\left\langle Y,  b\right\rangle_H
+\left\langle z,  g\right\rangle_{L_2(E_2;H)} +\left\langle Z,  \sigma
\right\rangle_{L_2(E_1;H)} +\left\langle \left\langle k,  \varphi\right\rangle
\right\rangle_{L_{\Pi  }^{2}(   H )} ,$$
and
$$
\left\langle \left\langle k,  \varphi\right\rangle \right\rangle_{L_{\Pi  }^{2}(   H )}  (t, \upsilon )
:=\left\langle \left\langle k,  \varphi (t, \upsilon ) \right\rangle \right\rangle_{L_{\Pi  }^{2}(   H )}  (t, \upsilon )
=\int_{\Theta}\left\langle k \left(  \rho\right)  ,  \varphi\left(
t,\upsilon ,\rho\right)  \right\rangle \Pi\left(  d\rho\right)  .
$$

Let us now set our hypotheses.

\medskip

(A1) (Monotonicity condition 1): ${\forall\, \upsilon=\left(
y,Y,z,Z,k\right), \upsilon'=\left( y',Y',z',Z'
,k'\right)\in \mathbb{H}^2}$, and $\forall \, t\in\left[  0,T\right]$,
\begin{align*}
\hspace{0.5cm} \left\langle A\left(  t,\upsilon\right)  -A\left(
t,\upsilon' \right)  ,\upsilon-\upsilon' \right\rangle   \leq
& -\theta_{1}\, \big( \left\vert y-y'  \right\vert_H^{2}+\left\Vert  z-z' \right\Vert^{2}_{L_2(E_2;H)}\big)  \\
& \hspace{-1cm} -\, \theta_{2}\, \big(  \left\vert  Y-Y' \right\vert_H^{2}+\left\Vert
 Z-Z'  \right\Vert^{2}_{L_2(E_1;H)}
 +||| k-k' |||^{2}_{L_{\Pi }^{2}( H )}\big)  .
\end{align*}

(A2) (Monotonicity condition 2):
 \begin{equation*}
 \left\langle h\left(  y\right)  -h\left(  y'\right)  ,  y-y'  \right\rangle_H  \geq\beta\left\vert   y-y' \right\vert_H^{2} \quad \quad \forall \; y,y'\in H.
\end{equation*}
Here, $\theta_{1},\theta_{2},$\ and $\beta$\ are
nonnegative constants such that ${\theta_{1}+\theta_{2}>0}$ and ${\theta_{2}+\beta>0.}$

(A3) For all $\upsilon\in\mathfrak{M}^{2}, \, A\left(
t,\upsilon\right)  $ is an $\mathcal{F}_{t}$-measurable $\mathbb{H}^2$-valued stochastic process defined
on $\Omega\times \left[  0,T\right]  $ with $A\left(  .,0\right)  \in
\mathfrak{M}^{2}$, and for all $y\in H ,h\left(  y\right) $ is $\mathcal{F}_{T}$-measurable with $h\left(  0\right)\in L^{2}\left(  \Omega,\mathcal{F}_{T},\mathbb{P}; H\right).$

(A4) (Lipschitz condition): $\exists \, c>0$ and $0<\gamma<1$  such that
\begin{eqnarray*}
&& \hspace{-0.5cm} \left\vert b\left(  t,y,Y,z,Z,k\right)  -b\left(  t,y',Y',z',Z',k'\right)  \right\vert^2 \leq c \, \big(  \left\vert y-y'\right\vert ^{2}+\left\vert Y-Y'\right\vert ^{2}
+\left\Vert z-z'\right\Vert^{2} \\  && \hspace{3.3in} +\,  \left\Vert Z-Z'\right\Vert^{2}+||| k-k'|||^{2}\big)  , \\
&&  \hspace{-0.5cm} \left\vert f\left(  t,y,Y,z,Z,k\right)  -f\left( t,y',Y',z',Z',k'\right)  \right\vert ^{2} \leq c \,  \big(   \left\vert y-y'\right\vert ^{2}+\left\vert Y-Y'\right\vert ^{2}+\left\Vert z-z'\right\Vert^{2} \\  && \hspace{3.3in} +\,  \left\Vert Z-Z'\right\Vert^{2} + ||| k-k'|||^{2}\big)   ,\\
&&   \hspace{-0.5cm}  \left\vert \sigma\left(  t,y,Y,z,Z,k\right)  -\sigma\left(
t,y',Y',z',Z',k'\right)  \right\vert ^{2}  \leq c \, \big(   \left\vert y-y'\right\vert ^{2}+\left\vert Y-Y'\right\vert ^{2}+\left\Vert Z-Z'\right\Vert^{2} \\  && \hspace{3.3in} + \, ||| k-k'||| ^{2} \big) +\frac{\gamma}{2} \left\Vert z-z'\right\Vert^{2}  , \\
&& \hspace{-0.5cm} \left\vert g\left(  t,y,Y,z,Z,k\right)  -g\left(
t,y',Y',z',Z',k'\right)  \right\vert ^{2} \leq c \left(  \left\vert y-y'\right\vert ^{2}+\left\vert Y-Y'\right\vert ^{2}+\left\Vert z-z'\right\Vert^{2}\right)
\\  \hspace{-0.5cm} && \hspace{3.3in}
+ \, \gamma\left(  \left\Vert Z-Z'\right\Vert^{2}+||| k-k'||| ^{2}\right)  ,\\
&&  \hspace{-0.5cm}   \int_{\Theta}\left\vert \varphi\left(  t,y,Y,z,Z,k,\rho
\right)  -\varphi\left( t,y',Y',z',Z',k'
,\rho\right)  \right\vert ^{2} \Pi(d\rho) \\ &&  \hspace{1.5cm}  \leq c \left(  \left\vert y-y'\right\vert ^{2}+\left\vert Y-Y'\right\vert ^{2}+ \left\Vert Z-Z'\right\Vert^{2}+||| k-k'|||^{2}\right)
+\frac{\gamma}{2} \left\Vert z-z'\right\Vert^{2}  ,\\
&&  \hspace{-0.5cm} \left\vert h\left(  y\right)  -h\left(
y'\right)  \right\vert \leq c\left\vert y-y'\right\vert ,
\end{eqnarray*}
for all $(y,Y,z,Z,k), \, (y',Y',z',Z',k') \in \mathbb{H}^2$ and all $t\in [0,T]$.
\begin{remark}\label{Remark 3.1}
It is useful to know that we can also replace hypotheses {\rm (A1)} and {\rm (A2)}
by the following ones and derive essentially the same results as in the following two theorems:

\smallskip

{\rm (A1)$'$} $\, \forall \,\upsilon=\left(  y,Y%
,z,Z,k\right)  ,\upsilon'=\left(y',Y',z',Z',k'\right)  \in \mathbb{H}^2,$ $\, \forall \,
t\in\left[  0,T\right]  $%
\begin{align*}
\hspace{0.7cm} \left\langle A\left(  t,\upsilon\right)  -A\left(
t,\upsilon'\right)  ,\upsilon-\upsilon'\right\rangle_H \geq & \; \theta_{1}\, \big(  \left\vert y-y'  \right\vert_H^{2}+\left\Vert   z-z'  \right\Vert^{2}_{L_2(E_2;H)} \big) \\
& \hspace{-1.3cm} +\, \theta_{2}\, \big(  \left\vert  Y-Y' \right\vert_H
^{2}+\left\Vert  Z-Z' \right\Vert^{2}_{L_2(E_1;H)}+|||  k-k'|||_{L_{\Pi  }^{2}(   H )}^{2}\big) ,
\end{align*}

\begin{align*}
\hspace{-0.6cm} {\rm (A2)'} \hspace{3cm}  \left\langle h\left(  y\right)  -h\left(  y'\right)  , y-y'  \right\rangle  & \leq-\beta\left\vert y-y' \right\vert ^{2} \quad \forall \; y,y'\in H.
\end{align*}
\end{remark}

\smallskip

In the following, we present our main results.
\begin{theorem}[Uniqueness]\label{Propo 3.4-1}
Under hypotheses {\rm(A1)--(A4)} (or {\rm(A1)$'$, (A2)$'$, (A3)--(A4)}; cf. 
Remark~\ref{Remark 3.1}), FBDSDEJ~(\ref{eq:3.1}) has at most one solution
$\left(y,Y,z,Z,k\right)$  in $\mathfrak{M}^{2}$.
\end{theorem}

For the proof of this theorem, one can closely follow the proof of Theorem~2.3 in \cite{Al-G-relaxed}, with minor modifications. Thus, we omit it in order to avoid increasing the size of this paper.

It is important to highlight that the inclusion of $\frac{\gamma}{2}$ in the Lipschitz conditions associated with the mappings $\sigma$ and $\varphi$ in {\rm (A4)} is necessary to establish the uniqueness of such solutions, as explained in \cite{Al-G-relaxed}. This aspect of uniqueness is also crucial for constructing unique solutions over small intervals using the method of time continuation, as we will clearly observe in the proof of Lemma~\ref{lem:final-lemma} below.

\begin{theorem}[Existence]\label{Propo 3.4}
Assume {\rm(A1)--(A4)} or {\rm((A1)$'$, (A2)$'$, (A3)--(A4))}. Then FBDSDEJ~(\ref{eq:3.1}) has a solution
$\left(y,Y,z,Z,k\right)$  in $\mathfrak{M}^{2}.$
\end{theorem}

We will concentrate exclusively on conditions {\rm(A1)--(A4)}, as the alternative case assuming {\rm((A1)$'$, (A2)$'$, (A3)--(A4))} can be handled in a similar manner without any difficulty.

To ensure clarity and organization, we divide the proof into two separate cases.

\bigskip

\noindent\textbf{Case~1:} Let $\theta_1 >0, \theta_2 \geq0$ and $\beta >0$. Consider the following family of FBDSDEJ parameterized by $\alpha\in\left[
0,1\right]  :$
\begin{eqnarray}\label{eq:3.3}
\left\{
\begin{array}{ll}
dy_{t}=[  \alpha \,  b\left(  t,\upsilon_{t}\right)  +\widetilde{b}%
_{0}(t)  ]  \, dt+[  \alpha \,  \sigma\left(
t,\upsilon_{t}\right)  +\widetilde{\sigma}_{0}(t)  ] \,
dW_{t}-z_{t}\, \overleftarrow{dB}_{t} \\ \hspace{2.5in}
+\int_{\Theta}[  \alpha \, \varphi\left(  t,\upsilon_{t},\rho\right)
+\varphi_{0}\left(  t,\rho\right)  ]  \widetilde{N}\left(
d\rho,dt\right), \\ \\
dY_{t}=[  \alpha \,  f\left(  t,\upsilon_{t}\right)  -\left(
1-\alpha\right)  \theta_{1} \,  y_{t}+\widetilde{f}_{0}(t)  ]\,
dt+Z_{t}\, dW_{t}\\  \hspace{1in}
+[  \alpha \,  g\left(  t,\upsilon_{t}\right)  -\left(
1-\alpha\right)  \theta_{1} \,  z_{t}+\widetilde{g}_{0}(t)  ]
\, \overleftarrow{dB}_{t}+\int_{\Theta}k_{t}\left(  \rho\right)  \widetilde
{N}(d\rho,dt)  ,
\\ \\
y_{0}=x,Y_{T}=\alpha \,  h\left(
y_{T}\right)  +\left(  1-\alpha\right)   y_{T}+\phi,
\end{array}
\right.
\end{eqnarray}
where $\left(
\widetilde{b}_{0},\widetilde{f}_{0},\widetilde{\sigma}_{0},\widetilde{g}_{0},\varphi
_{0}\right)  \in\mathfrak{M}^{2}$ and $\phi\in L^{2}\left(  \Omega
,\mathcal{F}_{T},\mathbb{P}; H \right)$ are given arbitrarily.

When $\alpha=1$, the existence of the solution of (\ref{eq:3.3}) implies clearly
that of (\ref{eq:3.2}) by taking $(\widetilde{b}_{0},\widetilde{f}_{0},\widetilde{\sigma}_{0},\widetilde{g}_{0},\varphi_{0})=(0,0,0,0,0)$, while when $\alpha=0$, (\ref{eq:3.3}) becomes a decoupled FBDSDEJ of the following shape:
\begin{eqnarray}\label{eq:3.4}
\hspace{-0.2cm}\left\{
\begin{array}{ll}%
dy_{t}=\widetilde{b}_{0}(t)\, dt+\widetilde{\sigma}_{0}(t)\, dW_{t}-z_{t}\, \overleftarrow{dB}_{t}+\int_{\Theta}\varphi_{0}\left(  t,\rho\right)\widetilde{N}\left(
d\rho,dt\right),\\ \\
dY_{t}=[ - \theta_{1} \,  y_{t} + \widetilde{f}_{0}(t) ] \, dt + Z_{t}\, dW_{t}+[ -\theta_{1} \,  z_{t}+\widetilde{g}_{0}(t)  ]
\, \overleftarrow{dB}_{t}  +\int_{\Theta}k_{t}\left(  \rho\right)  \widetilde{N}(d\rho,dt),
\\ \\
y_{0}=x,Y_{T}=  y_{T}+\phi, \; 0\leq t\leq T.
\end{array}
\right.
\end{eqnarray}

Next, by making use of (\ref{eq:2.3}) (with $W$ replacing $\breve{B}$) and 
(\ref{eq:2.6}), one can rewrite
\( \int_{0}^{t}\widetilde{\sigma}_{0}(s)\, dW_{s}$ as $-\int_{T-t}^{T}\breve{\widetilde{\sigma}}_{0}(s)\, \overleftarrow{d\breve{W}}_{s},\) where 
$\breve{\widetilde{\sigma}}_{0}(s)=\widetilde{\sigma}_{0}(T-s)$, i.e., as a backward It\^{o} integral. Furthermore, \(\int_{0}^{t}z_{s}\, \overleftarrow{dB}_{s}=-\int_{T-t}^{T}\breve{z}_{s}\, d\breve{B}_{s},\)
with $\breve{Z}_{s}=Z_{T-s}$, yielding a forward It\^{o} integral. The same applies to the Lebesgue integrals (with respect to $dt$) and the integrals with respect to $\widetilde{N}$. This business helps to rewrite the first (forward) equation of (\ref{eq:3.4}) when replacing $t$ there by $T-t$ as a BDSDE  with jumps (BDSDEJ) of the form:
\begin{equation}\label{eq:BDSDEJ}
\breve{y}_t  =\breve{y}_T - \int_{t}^{T}\breve{\widetilde{b}}_{0}(s)\, \overleftarrow{ds} - \int_{t}^{T}\breve{\widetilde{\sigma}}_{0}(s)\, \overleftarrow{d\breve{W}}_{s} - \int_{t}^{T}\int_{\Theta}\breve{\varphi}_{0}\left(  s,\rho\right)\overleftarrow{\breve{\widetilde{N}}}(d\rho,ds) +\int_{t}^{T}\breve{z}_{s}\, d\breve{B}_{s},
\end{equation}
$0\leq  t\leq T,$ where $\breve{y}_t=y_{T-t},$ (and so $\breve{y}_{T}=y_{0}),\breve{\widetilde{b}}_{0}(s)=\widetilde{b}(T-s),\breve{\varphi}_{0}\left(  s,\rho\right)=\varphi(T-s,\rho)$ and
$\overleftarrow{ds}= - ds$.
Consequently, under (A3) it follows that (\ref{eq:BDSDEJ}) attains a unique solution $(\breve{y},\breve{z})$ in $\mathcal{M}^{2}\left( 0,T;H\right)
\times\mathcal{M}^{2}\left(0,T;L_2(E_2;H)\right)$. Namely,
\[
\breve{y}_t= \mathbb{E} \Big[ x - \int_t^T \breve{\widetilde{b}}_{0}(s)\, \overleftarrow{ds}-\int_{t}^{T}\breve{\widetilde{\sigma}}_{0}(s)\, \overleftarrow{d\breve{W}}_{s}-\int_{t}^{T}\int_{\Theta}\breve{\varphi}_{0}\left(  s,\rho\right)\overleftarrow{\breve{\widetilde{N}}}(d\rho,ds) \, \big{|} \, \mathcal{F}_t \Big]
\]
and $\breve{z}$ is got via the generalized martingale representation theorem (\cite[Theorem~1]{AG}) by the following formula:
\begin{eqnarray*}
M(t) = M(0) + \int_{0}^{t}(-\breve{z}_{s})\, d\breve{B}_{s}  = M(0) - \int_{0}^{t} \breve{z}_{s} \, d\breve{B}_{s} ,
\end{eqnarray*}
where
$$M(t)=\mathbb{E} \Big[ x-\int_0^T \breve{\widetilde{b}}_{0}(s)\, \overleftarrow{ds}-\int_{0}^{T}\breve{\widetilde{\sigma}}_{0}(s)\, \overleftarrow{d\breve{W}}_{s}-\int_{0}^{T}\int_{\Theta}\breve{\varphi}_{0}(s,\rho)\overleftarrow{\breve{\widetilde{N}}}(d\rho,ds) \, \big{|} \, \mathcal{F}_t \Big] .$$
Therefore, $\{(y_{s},z_{s}):=(\breve{y}_{T-s},\breve{z}_{T-s}),0\leq t\leq T\}$ is the unique solution of the forward equation of (\ref{eq:3.4}).

Now, by substituting $(y,z)$ into the second (backward) equation of (\ref{eq:3.4}), it remains a BDSDEJ of type (\ref{eq:BDSDEJ}), and thus it admits a unique solution 
$(Y,Z,k)$. Specifically,
\[
Y_t= \mathbb{E} \Big[ y_{T}+\phi - \int_t^T ( \theta_{1} \,  y_{s} + \widetilde{f}_{0}(s) ) \, ds - \int_t^T ( -\theta_{1} \,  z_{s}+\widetilde{g}_{0}(s) ) \, \overleftarrow{dB}_{s} \, \big{|} \, \mathcal{F}_t \Big],
\]
for each $t\in [0,T]$, and the two processes $(Z,k)$ can be obtained again via the generalized martingale representation theorem (\cite[Theorem~1]{AG}) as follows:
\begin{eqnarray*}
P(t)  = P(0) + \int_0^t Z_{s}\, dW_{s} + \int_0^t \int_{\Theta}k_{s}\left(  \rho\right)  \widetilde{N}(d\rho,ds) ,
\end{eqnarray*}
where, $0 \leq t \leq T$,
 $$P(t)= \mathbb{E} \Big[ y_{T}+\phi - \int_0^T ( \theta_{1} \,  y_{s} + \widetilde{f}_{0}(s) ) \, ds - \int_0^T ( -\theta_{1} \,  z_{s}+\widetilde{g}_{0}(s)  ) \, \overleftarrow{dB}_{s} \, \big{|} \, \mathcal{F}_t \Big] .$$
We have then derived a unique solution $\upsilon=(y,Y,z,Z,k)$ of (\ref{eq:3.4}) in $\mathfrak{M}^{2}.$

\medskip

Our next step is to deal with system~(\ref{eq:3.1}). We shall use the method of continuation. We start with a priori lemma.
\begin{lemma}\label{Lemma: 3.5}
Assume that $\theta_1 >0, \theta_2 \geq 0$ and $\beta >0$. Under hypotheses {\rm(A1)--(A4)} there exists a positive constant $\delta_{0}$\ that is independent of $\alpha_0$ such
that if, a priori, for each $\phi\in L^{2}\left(  \Omega,\mathcal{F}_{T},\mathbb{P};H\right)  $\ and $\left(  \widetilde{b}_{0},\widetilde{f}_{0},\widetilde{g}_{0},\widetilde{\sigma}_{0},\varphi_{0}\right)  \in\mathfrak{M}^{2}$,
 FBDSDEJ~(\ref{eq:3.3}) is uniquely solvable for some fixed $\alpha_{0}\in[  0,1)  $, then
for each $\alpha\in\left[  \alpha_{0},\alpha_{0}+\delta_{0}\right] ,\phi\in
L^{2}\left(  \Omega,\mathcal{F}_{T},\mathbb{P};H \right)  $\ and $\left(
\widetilde{b}_{0},\widetilde{f}_{0},\widetilde{\sigma}_{0},\widetilde{g}_{0},\varphi
_{0}\right)  \in\mathfrak{M}^{2}$, this system~(\ref{eq:3.3}) is uniquely solvable.
\end{lemma}
\begin{proof}
Let us assume that, for any $\xi \in L^{2}\left( \Omega ,\mathcal{F}_{T},\mathbb{P}
,H\right) ,$ $\left(  \widetilde{b}_{0},\widetilde{f}_{0},\widetilde{g}_{0},\widetilde{\sigma}_{0},\varphi_{0}\right)  \in\mathfrak{M}^{2}$, the
FBDSDEJ~(\ref{eq:3.2})
has a unique solution for a (fixed) constant $\alpha =\alpha _{0}\in \left[ 0,1\right)$. Then, for each element $\bar{v}%
=\left( \bar{y},\bar{Y},\bar{z},\bar{Z}, \bar{k} \right)$  of  $\mathfrak{M}^{2},$  there exists
 $v=\left( y,Y,z,Z,k\right)\in \mathfrak{M}^{2}$ solving uniquely the FBDSDEJ:
\begin{eqnarray}\label{eq:new-number}
\left\{
\begin{array}{ll}%
\hspace{-0.25cm} dy_{t}=\left[  \alpha_{0}\, b\left(  t,\upsilon_{t}\right)
+\delta \,  b\left(  t,\overline{\upsilon}_{t}\right)  +\widetilde{b}_{0}\left(
t\right)  \right]  dt-z_{t}\, \overleftarrow{dB}_{t}\\ \hspace{1.5cm}
+\left[  \alpha_{0}\, \sigma\left(  t,\upsilon_{t}\right)  +\delta \,
\sigma\left(  t,\overline{\upsilon}_{t}\right)  +\widetilde{\sigma}_{0}\left(
t\right)  \right]  dW_{t}\\ \hspace{1.5cm}
+\int_{\Theta}\left[  \alpha_{0}\, \varphi\left(  t,\upsilon_{t}%
,\rho\right)  +\delta \, \varphi\left(  t,\overline{\upsilon}_{t},\rho\right)
+\varphi_{0}\left(  t,\rho\right)  \right]  \widetilde{N}\left(  d\rho,dt\right),\\ \\
\hspace{-0.25cm} dY_{t}=\left[  \alpha_{0}\, f\left(  t,\upsilon_{t}\right)  -\left(
1-\alpha_{0}\right)  \theta_{1} \, y_{t}+\delta  \left(  f\left(
t,\overline{\upsilon}_{t}\right)  +\theta_{1}\, \overline{y}_{t}\right)  +\widetilde{f}%
_{0}\left(  t\right)  \right]dt\\
\hspace{1.5cm} +\left[  \alpha_{0}\, g\left(  t,\upsilon_{t}\right)  -\left(
1-\alpha_{0}\right)  \theta_{1}\, z_{t}+\delta \left(  g\left(
t,\overline{\upsilon}_{t}\right)  +\theta_{1} \, \overline{z}_{t}\right)  +\widetilde{g}%
_{0}\left(  t\right)  \right]  \overleftarrow{dB}_{t} \\ \hspace{7cm}
+\, Z_{t}\, dW_{t}+\int_{\Theta}k_{t}\left(  \rho\right)  \widetilde{N}\left(
d\rho,dt\right) ,
\\ \\
y_{0}=x ,
Y_{T}=\alpha_{0}\, h\left(  y_{T}\right)  +\left(  1-\alpha_{0}\right)
y_{T}+\delta \left(  h\left(  \overline{y}_{T}\right)  - \overline{y}
_{T}\right)  +\phi.
\end{array}
\right.
\end{eqnarray}

\bigskip

Define the mapping $I_{\alpha_{0}+\delta}$ by
\begin{eqnarray*}
I_{\alpha_{0}+\delta}\left(
\overline{\upsilon}_{\cdot},\overline{y}_{T}\right) := \left(  \upsilon_{\cdot},y_{T}\right) :\mathfrak{M}^{2}\times
L^{2}\left(  \Omega,\mathcal{F}_{T},\mathbb{P};H\right) \rightarrow\mathfrak{M}^{2}\times L^{2}\left(  \Omega,\mathcal{F}%
_{T},\mathbb{P};H\right) .
\end{eqnarray*}
We shall show that the mapping $I_{\alpha_{0}+\delta}$ is a contraction for some small $\delta >0$.
Let
\(
\overline{\upsilon}_{\cdot} :=\left(  \overline{y}_{\cdot},\overline{Y}
_{\cdot},\overline{z}_{\cdot},\overline{Z}_{\cdot},\overline{k}_{\cdot}\right), \overline{\upsilon}_{\cdot}^{\prime} :=\left(  \overline{y}_{\cdot}^{\prime},\overline{Y}
_{\cdot}^{\prime},\overline{z}_{\cdot}^{\prime},\overline{Z}_{\cdot}^{\prime},\overline{k}_{\cdot}^{\prime
}\right)$ be two elements of $\mathfrak{M}^{2}$.
Then let
$$I_{\alpha_{0}+\delta
}\left(  \overline{\upsilon}_{\cdot}, \overline{y}_{T}\right):=\left({\upsilon}_{\cdot},{y}_{T}\right) =
(\left(  y_{\cdot},Y_{\cdot},z_{\cdot},Z_{\cdot},k_{\cdot}\right), y_{T})$$
and
$$I_{\alpha_{0}+\delta
}\left(  \overline{\upsilon}^{\prime}_{\cdot}, \overline{y}^{\prime}_{T} \right):= (\upsilon_{\cdot}^{\prime}, y^{\prime}_{T})=(\left(  y_{\cdot}^{\prime},Y_{\cdot}^{\prime},z_{\cdot}^{\prime},Z_{\cdot}^{\prime},k_{\cdot}^{\prime}\right), y^{\prime}_{T}) .$$
To simplify the notation, let us denote
$$
\widehat{\overline{\upsilon}}_{t}=\overline{\upsilon}_{t}-\overline{\upsilon}_{t}^{\prime}=\left(  \widehat{\overline{y}}_{t},\widehat{\overline{Y}}_{t},
\widehat{\overline{z}}_{t},\widehat{\overline{Z}}_{t},\widehat{\overline{k}}_{t}\right),
\; \; \widehat{\upsilon}_{t}=\upsilon_{t}-\upsilon_{t}^{\prime}=\left(  \widehat{y}_{t},\widehat{Y}_{t},\widehat{z}_{t},\widehat{Z}_{t},\widehat{k}_{t}\right),$$
$$\widehat{\pi}\left(  t,\overline{\upsilon}_{t}\right):=\pi
\left(  t,\overline{\upsilon}_{t}\right) - \pi \left(  t,\overline{\upsilon}
_{t}^{\prime}\right),
\;\;
\widehat{\pi}\left(  t,\upsilon_{t}\right)
:=\pi\left(  t,\upsilon_{t}\right)  -\pi \left(  t,\upsilon
_{t}^{\prime}\right)$$ for ${\pi  :=b,\sigma,f,g}$, and also $$\widehat{h}(\overline{y}_T):= h(\overline{y}_T)-h(\overline{y}_{t}^{\prime}) , \;\; \widehat{h}(y_T):= h(y_T)-h(y_{t}^{\prime})$$ and
 $$\widehat{\varphi}\left(  t, \overline{\upsilon}_{t},\cdot\right):=\varphi
\left(  t, \overline{\upsilon}_{t},\cdot\right)  -\varphi\left(  t, \overline{\upsilon}_{t}^{\prime},\cdot\right) , \;\;
\widehat{\varphi}\left(  t, \upsilon_{t},\cdot\right):=\varphi
\left(  t, \upsilon_{t},\cdot\right)  -\varphi\left(  t, \upsilon_{t}^{\prime},\cdot\right) .$$

With the help of hypothesis (A3), by applying It\^{o}'s formula to $\langle \widehat{y},\widehat{Y}\rangle
$\ on $\left[  0,T\right]$ and taking the expectation, we find that
\begin{eqnarray}\label{eq:3.7}
&& \hspace{-1.25cm} \mathbb{E} \left[  \left\langle \widehat{y}_{T},\alpha_{0}\, \widehat{h}\left(
y_{T}\right) +\left(  1-\alpha_{0}\right)  \widehat{y}_{T}\right\rangle \right]  +\left(  1-\alpha_{0}\right)  \theta_{1}\, \mathbb{E} \left[  \int_{0}^{T}\left(
\left\vert \widehat{y}_{t}\right\vert ^{2}+\left\Vert \widehat{z}_{t}\right\Vert
^{2}\right)  dt\right]\nonumber \\
&& \hspace{4.8cm}-\, \alpha_{0} \, \mathbb{E} \left[  \int_{0}^{T}\left\langle A\left(
t,\upsilon_{t}\right)  -A\left(  t,\upsilon_{t}^{\prime}\right)
,\widehat{\upsilon}_{t}\right\rangle dt\right]
\nonumber \\
&&  =-\, \delta\, \mathbb{E} \left[  \left\langle \widehat{y}_{T},\widehat{h}\left(
\overline{y}_{T}\right)  -\widehat{\overline{y}}_{T}\right\rangle \right]+ \delta\, \mathbb{E} \left[  \int_{0}^{T}\left\langle \widehat{y}_{t}%
,\widehat{f}\left(  t,\overline{\upsilon}_{t}\right)  \right\rangle
dt\right] \nonumber \\
&& \hspace{1cm} + \, \delta\, \theta_{1}\,  \mathbb{E} \left[  \int_{0}^{T}\left\langle
\widehat{y}_{t},\widehat{\overline{y}}_{t}\right\rangle dt\right]+ \delta\, \mathbb{E} \left[  \int_{0}^{T}\left\langle \widehat{Y}_{t},\widehat
{b}\left(  t,\overline{\upsilon}_{t}\right)  \right\rangle dt\right] \nonumber \\
&& \hspace{1cm}
+\, \delta\, \mathbb{E} \left[  \int_{0}^{T}\left\langle \widehat{\sigma%
}\left(  t,\overline{\upsilon}_{t}\right)  ,\widehat{Z}_{t}\right\rangle dt\right]+ \delta\, \mathbb{E} \left[  \int_{0}^{T}\left\langle \widehat{z}_{t}%
,\widehat{g}\left(  t,\overline{\upsilon}_{t}\right)  \right\rangle
dt\right] \nonumber\\
&& \hspace{1cm} + \, \delta\, \theta_{1}\, \mathbb{E} \left[  \int_{0}^{T}\left\langle
\widehat{z}_{t},\widehat{\overline{z}}_{t}\right\rangle dt\right]+ \delta\, \mathbb{E} \left[  \int_{0}^{T}\left\langle \left\langle
\widehat{\varphi}\left(  t,\overline{\upsilon}_{t},\cdot\right)  ,\widehat{k}_{t
}\right\rangle \right\rangle dt\right].
\end{eqnarray}
Then, by making use of (A2) and (A4), it follows that
\begin{eqnarray*}
&& \hspace{-1.5cm} \left(  1-\alpha_{0}+\alpha_{0}\, \beta\right)  \, \mathbb{E} \left[
\left\vert \widehat{y}_{T}\right\vert ^{2}\right]  +\theta_{1}\, \mathbb{E} \left[ \int
_{0}^{T}\left(  \left\vert \widehat{y}_{t}\right\vert ^{2}+\left\Vert \widehat{z}_{t}\right\Vert^{2}\right)  dt \right] \nonumber
\\
&& \hspace{1cm}+ \, \alpha_{0}\, \theta_{2}\, \mathbb{E} \left[  \int_{0}^{T}\left(  \vert
\widehat{Y}_{t}\vert ^{2}+\Vert \widehat{Z}_{t}\Vert
^{2}+||| \widehat{k}_{t}||| ^{2}\right)  dt\right]
\nonumber \\ &&
\leq  \delta  \left(  \mathbb{E} \left[  \left\vert \widehat{y}_{T}\right\vert
^{2}\right]  + C\, \mathbb{E} \left[  \vert \widehat{\overline{y}}_{T}\vert
^{2}\right]  \right) + \, \delta \, C\, \mathbb{E} \left[ \int_{0}^{T}\left( \Vert \widehat{\upsilon}%
_{t}\Vert^{2}+\Vert \widehat{\overline{\upsilon}}_{t}\Vert^{2}\right)  dt \right],
\end{eqnarray*}
for some constant $C>0$, which can vary from place to place and depends only on the constants $c,\gamma,\theta_{1},\theta_{2}, \beta$, and $T$.
Since $$1-\alpha_{0}+\alpha_{0}\, \beta\geq \min \{\beta , 1\}=:\tilde{\beta}$$
and $\tilde{\beta} > 0$, this inequality implies
\begin{eqnarray*}
&& \hspace{-1.5cm} \tilde{\beta} \, \,  \mathbb{E} \left[
\left\vert \widehat{y}_{T}\right\vert ^{2}\right]  + \theta_1\,  \mathbb{E} \left[
\int_{0}^{T}\left(  \vert \widehat{y}_{s}\vert ^{2}+\left\Vert
\widehat{z}_{s}\right\Vert^{2}\right)  ds\right]
\nonumber\\
&&  \leq\delta \, C \left(  \mathbb{E} \left[  \left\vert \widehat{y}%
_{T}\right\vert ^{2}\right]  +\mathbb{E} \left[  \vert \widehat{\overline{y}
}_{T}\vert ^{2}\right]     +\, \mathbb{E} \left[  \int_{0}^{T}\left(  \left\Vert \widehat{\upsilon
}_{s}\right\Vert^{2}+\Vert \widehat{\overline{\upsilon}}%
_{s}\Vert^{2}\right)  ds\right]\right),
\end{eqnarray*}
after ignoring the term $\alpha_{0}\, \theta_{2}\, \mathbb{E} \Big[  \int_{0}^{T}\big(  \vert
\widehat{Y}_{t}\vert ^{2}+\Vert \widehat{Z}_{t}\Vert
^{2}+||| \widehat{k}_{t}||| ^{2}\big)  dt \Big]$ in which $\alpha_0$ appears because $\delta$ we are looking should be independent of $\alpha_0$.
Now, letting ${\lambda =\min\{  \tilde{\beta},\theta_{1}\}}$, we realize that
$0< \lambda  < 1$ and
\begin{eqnarray}\label{eq:3.11}
&& \hspace{-1.5cm} \mathbb{E} \left[  \left\vert \widehat{y}_{T}\right\vert
^{2}\right]  +\mathbb{E} \left[  \int_{0}^{T}\left(  \left\vert \widehat{y}%
_{s}\right\vert ^{2}+\left\Vert \widehat{z}_{s}\right\Vert^{2}\right)
ds \right] \nonumber
\\
&& \leq \frac{\delta \, C}{\lambda } \left(  \mathbb{E} \left[  \left\vert \widehat{y}%
_{T}\right\vert ^{2}\right]  +\mathbb{E} \left[  \vert \widehat{\overline{y}
}_{T}\vert ^{2}\right]    + \mathbb{E} \left[  \int_{0}^{T}\left(  \left\Vert \widehat{\upsilon
}_{s}\right\Vert^{2}+\Vert \widehat{\overline{\upsilon}}
_{s}\Vert^{2}\right)  ds\right]\right).
\end{eqnarray}

Next, we want to find a similar estimate for $\mathbb{E} \left[  \int_{0}^{T}\left(  \vert \widehat{Y}_{s}\vert
^{2}+\Vert \widehat{Z}_{s}\Vert^{2}+||| \widehat{k}_{s}||| ^{2}\right)  ds\right]$. We apply first It\^{o}'s formula to $\vert \widehat{Y}_{s}\vert ^{2}$ and use Cauchy-Schwarz inequality to obtain
\begin{eqnarray}\label{eq:3.12}
&& \hspace{-0.75cm}\mathbb{E} \left[  \vert \widehat{Y}_{t}\vert ^{2}\right]
+\mathbb{E} \left[  \int_{t}^{T}\Vert \widehat{Z}_{s}\Vert
^{2}ds\right]  +\mathbb{E} \left[  \int_{t}^{T}||| \widehat
{k}_{s}||| ^{2} \, ds\right] \nonumber \\ && \hspace{-0.5cm}
\leq 4\, \mathbb{E} \left[  \alpha_{0}^{2}\, \left\vert \widehat{h}\left(
y_{T}\right)  \right\vert ^{2}+\left(  1-\alpha_{0}\right)  ^{2}\left\vert
\widehat{y}_{T}\right\vert ^{2}+\delta^{2} \, \left\vert \widehat{h}\left(  \overline{y}_{T}\right)  \right\vert ^{2}+\delta^{2}\, \vert \widehat{\overline{y}}%
_{T}\vert ^{2}\right]
\nonumber \\ && \hspace{-0.5cm}
+ \, 2\, \mathbb{E} \left[  \int_{t}^{T}\vert \widehat{Y}_{s}\vert
\cdot\left(  \alpha_{0}\, \left\vert \widehat{f}\left(  s,\upsilon
_{s}\right)  \right\vert +\left(  1-\alpha_{0}\right)  \theta_{1}\, \left\vert
\widehat{y}_{s}\right\vert + \delta \left\vert \widehat{f}\left(
s,\overline{\upsilon}_{s}\right)  \right\vert + \delta \, \theta_{1}\, \vert
\widehat{\overline{y}}_{s}\vert \right)  ds\right]  \nonumber\\
&& \hspace{3cm} +\, \mathbb{E} \left[  \int_{t}^{T}\left(\frac{1+\gamma}{2\gamma}\right)\alpha_{0}^{2}\left\vert \widehat{g}\left(  s,\upsilon_{s}\right)  \right\vert
^{2}ds\right]\nonumber \\ \hspace{-0.5cm}
&& + \, 3\, \mathbb{E} \left[  \int_{t}^{T}\left(  \frac{1+\gamma}{1-\gamma
}\right)  \left(  \left(  1-\alpha_{0}\right)  ^{2}\theta_{1}^{2}\left\Vert \widehat
{z}_{s}\right\Vert^{2} + \delta^{2}\left\vert \widehat{g}\left(
s,\overline{\upsilon}_{s}\right)  \right\vert ^{2}+\delta^{2}\, \theta_{1}^{2}\, \Vert \widehat{\overline{z}}_{s}\Vert^{2}\right)  ds\right]\nonumber \\
&& \hspace{-0.5cm} =: I_{1}+I_{2}(t)+I_{3}(t)+I_{4}(t),
\end{eqnarray}
where $I_i$ represents the $i$-th quantity in the right hand side of this inequality for $1\leq i \leq 4$.

We shall  now proceed to estimate each of these quantities. By applying (A3), (A4), and the Cauchy-Schwarz inequality, we obtain the following four inequalities:
\begin{equation}\label{eq:3.13}
 I_{1} \leq C\, \mathbb{E} \left[  \left\vert \widehat{y}_{T}\right\vert ^{2}
+\delta \, \vert \widehat{\overline{y}}_{T}\vert ^{2}\right],
\end{equation}
\begin{eqnarray}\label{eq:3.14}
&& \hspace{-0.5cm} I_{2} (t)  \leq C\, \mathbb{E} \left[  \int_{t}^{T}\vert \widehat{Y}_{s}\vert ^{2}\, ds\right]  +C\, \mathbb{E} \left[  \int_{t}^{T}\left\vert
\widehat{y}_{s}\right\vert ^{2}ds\right]  + C\, \delta\, \mathbb{E} \left[  \int_{t}^{T}\vert \widehat{\overline{y}}_{s}\vert ^{2}ds\right]  \nonumber \\ &&
\hspace{1cm} + \, C \, \mathbb{E} \left[  \int_{t}^{T}\left\Vert \widehat{z}_{s}\right\Vert^{2}ds\right]\nonumber + \; \left(\frac{1-\gamma}{8}\right)\, \mathbb{E} \left[  \int_{t}^{T}\left( \Vert \widehat
{Z}_{s}\Vert^{2}+||| \widehat{k}_{s}||| ^{2}\right)  ds\right]  \\ &&
\hspace{2.7in}  + \, C\, \delta\, \mathbb{E} \left[  \int_{t}^{T}\Vert \widehat{\overline{\upsilon}}_{s}\Vert^{2}ds\right],
\end{eqnarray}
\begin{eqnarray}\label{eq:3.15}
&& \hspace{-1.7cm} I_{3}(t) \leq C\, \mathbb{E} \left[  \int_{t}^{T}\left(  \left\vert \widehat{y}%
_{s}\right\vert ^{2}+\vert \widehat{Y}_{s}\vert ^{2}+\left\Vert \widehat
{z}_{s}\right\Vert^{2}\right)  ds\right] \nonumber  \\ &&
\hspace{2cm} +\left(
\frac{1+\gamma}{2}\right)  \alpha_{0}^{2}\;\mathbb{E} \left[ \int_{t}^{T} \left(  \Vert \widehat{Z}
_{s}\Vert^{2}+||| \widehat{k}
_{s}|||^{2}\right)  ds\right],
\end{eqnarray}
and
\begin{eqnarray}\label{eq:3.16}
&& \hspace{-1cm}  I_{4}(t)\leq C \, \mathbb{E} \left[  \int_{t}^{T}\left(  \left\Vert \widehat{z}%
_{s}\right\Vert^{2}+\delta \, \vert \widehat{\overline{y}}_{s}\vert
^{2}+\delta \, \Vert \widehat{\overline{z}}_{s}\Vert^{2}+\delta \, \vert
\widehat{\overline{Y}}_{s}\vert ^{2}\right)  ds\right]\nonumber\\
&& \hspace{2.5cm} +\, \delta \, C \left(  \frac{1+\gamma}{1-\gamma}\right)  \gamma \  \mathbb{E} \left[
\int_{t}^{T}\left( \Vert \widehat{\overline{Z}}_{s}\Vert^{2}
+||| \widehat{\overline{k}}_{s}||| ^{2}\right)  ds\right]  .
\end{eqnarray}
Then, we substitute (\ref{eq:3.13})-(\ref{eq:3.16}) in (\ref{eq:3.12}), noting that $\alpha_{0}<1$, to deduce that there exists a universal constant $C>0$, which is evidently independent of $\alpha_0$, such that
\begin{eqnarray}\label{eq:3.17}
&& \hspace{-1.5cm}\mathbb{E} \left[  \vert \widehat{Y}_{t}\vert ^{2}\right]  +\left(
1-\left(\frac{5+3\gamma}{8}\right)\right)  \mathbb{E} \left[  \int_{t}^{T}\left(  \Vert
\widehat{Z}_{s}\Vert^{2}+||| \widehat{k}_{s}||| ^{2}\right)  ds\right] \nonumber\\ &&
\leq C\, \mathbb{E} \left[  \int_{t}^{T}\vert \widehat{Y}_{s}\vert
^{2}ds\right]  +C\left(  \mathbb{E} \left[  \left\vert \widehat{y}_{T}\right\vert
^{2}\right]  +\delta\, \mathbb{E} \left[  \vert \widehat{\overline{y}}
_{T}\vert ^{2}\right]  \right) \nonumber\\ &&
 \hspace{3cm}+ \, C\, \mathbb{E} \left[  \int_{t}^{T}\left(
\left\vert \widehat{y}_{s}\right\vert ^{2}+\left\Vert \widehat{z}_{s}\right\Vert
^{2}+\delta  \Vert \widehat{\overline{\upsilon}}_{s}\Vert ^{2}\right)  ds\right] .
\end{eqnarray}
Since $0< \gamma <1$, it follows from applying Gronwall's inequality to (\ref{eq:3.17})  that
\begin{eqnarray}\label{eq:3.18}
&& \hspace{-1.5cm}\mathbb{E} \left[  \vert \widehat{Y}_{t}\vert ^{2}\right]  \leq
Ce^{T-t}\left(  C\left(  \mathbb{E} \left[  \left\vert \widehat{y}_{T}\right\vert
^{2}\right]  +\delta\, \mathbb{E} \left[  \vert \widehat{\overline{y}}
_{T}\vert ^{2}\right]  \right) \right.\nonumber \\ &&
\hspace{2.5cm}
\left. +\, C\, \mathbb{E} \left[  \int_{t}^{T}\left(
\left\vert \widehat{y}_{s}\right\vert ^{2}+\left\Vert \widehat{z}_{s}\right\Vert
^{2}+\delta \, \Vert \widehat{\overline{\upsilon}}_{s}\Vert ^{2}\right)  ds\right]  \right)  ,
\end{eqnarray}
for all $t\in [0,T]$.
Consequently,  (\ref{eq:3.17}) and (\ref{eq:3.18}) imply
\begin{eqnarray}\label{eq:3.22}
&& \hspace{-2cm} \mathbb{E} \left[  \int_{0}^{T}\left(  \vert \widehat{Y}_{s}\vert
^{2}+\Vert \widehat{Z}_{s}\Vert^{2}+||| \widehat{k}_{s}||| ^{2}\right)  ds\right]  \leq C\left(  \mathbb{E} \left[  \left\vert \widehat{y}
_{T}\right\vert ^{2}\right]  +\delta\, \mathbb{E} \left[  \vert \widehat
{\overline{y}}_{T}\vert ^{2}\right]  \right)  \nonumber \\
&& \hspace{3.7cm} +\, C\, \mathbb{E} \left[  \int_{0}^{T}\left(  \left\vert \widehat{y}_{s}\right\vert
^{2}+\left\Vert \widehat{z}_{s}\right\Vert^{2}+\delta \, \Vert \widehat
{\overline{\upsilon}}_{s}\Vert^{2}\right)  ds\right].
\end{eqnarray}

Inequality (\ref{eq:3.22}), along with (\ref{eq:3.11}), is crucial for establishing the proof of the lemma. Indeed, by combining these two inequalities, we obtain the following result:
\begin{eqnarray*}
&& \hspace{-1.5cm} \mathbb{E} \left[  \int_{0}^{T}\left\Vert \widehat{\upsilon}_{s}\right\Vert
^{2}ds\right]  +\mathbb{E} \left[  \left\vert
\widehat{y}_{T}\right\vert ^{2}\right] \nonumber
\\
&& \leq \frac{\delta \, C}{\lambda } \left(  \mathbb{E} \left[  \vert
\widehat{y}_{T}\vert ^{2}\right]  +\mathbb{E} \left[  \vert \widehat{\overline{y}
}_{T} \vert ^{2}\right]    + \mathbb{E} \left[  \int_{0}^{T}\left( \Vert \widehat{\upsilon}_{s}\Vert^{2}+\Vert
\widehat{\overline{\upsilon}}_{s}\Vert^{2}\right)  ds\right]\right)  \nonumber \\
&& \hspace{3.5cm} + \,
C\,\mathbb{E} \left[  \vert \widehat{y}%
_{T}\vert ^{2}\right]  + C\, \mathbb{E} \left[  \int_{0}^{T}\left(  \left\vert \widehat{y}_{s}\right\vert
^{2}+\left\Vert \widehat{z}_{s}\right\Vert^{2}\right)  ds\right].
\end{eqnarray*}
Hence, applying (\ref{eq:3.11}) once more to the last two terms of this latter inequality gives
\begin{eqnarray*}
&& \hspace{-1.5cm} \mathbb{E} \left[  \int_{0}^{T}\left\Vert \widehat{\upsilon}_{s}\right\Vert
^{2}ds\right]  +\mathbb{E} \left[  \left\vert
\widehat{y}_{T}\right\vert ^{2}\right]  \nonumber
\\
&& \leq \frac{\delta \, C}{\lambda } \left(  \mathbb{E} \left[
 \left\vert \widehat{y}_{T}\right\vert ^{2}\right]  +\mathbb{E} \left[  \vert \widehat{\overline{y}}_{T}\vert ^{2}\right] + \mathbb{E} \left[  \int_{0}^{T}\left(  \left\Vert \widehat{\upsilon}_{s}\right\Vert^{2}+\Vert \widehat{\overline{\upsilon}}
_{s}\Vert^{2}\right)  ds\right]\right) .
\end{eqnarray*}
Now take $\delta$ so that $0< \delta\leq \delta_0:=\frac{\lambda }{3C}$  to see that
\begin{equation*}
\mathbb{E} \left[  \int_{0}^{T}\left\Vert \widehat{\upsilon}_{s}\right\Vert
^{2}ds\right]  +\mathbb{E} \left[  \left\vert
\widehat{y}_{T}\right\vert ^{2}\right] \leq\frac{1}{2}\left(  \mathbb{E} \left[  \int_{0}^{T}\Vert \widehat
{\overline{\upsilon}}_{s}\Vert^{2}\, ds\right]  +\mathbb{E} 
\left[  \vert \widehat{\overline{y}}_{T}\vert ^{2}\right]
\right)  .
\end{equation*}

As a result, the mapping $I_{\alpha_{0}+\delta}$\ is a contraction for any fixed $\delta$ in $\left[ 0,\delta_{0}\right]$ on the Banach space $\mathfrak{M}^{2}\times L^{2}(\Omega,\mathcal{F}_{T},\mathbb{P};H)$.
Therefore, $I_{\alpha_{0}+\delta}$ attains a unique fixed point $\upsilon=\left(
y,Y,z,Z,k\right)$  in $\mathfrak{M}^{2},$ providing the unique solution to FBDSDEJ~(\ref{eq:3.3}) for $\alpha=\alpha_{0}+\delta,$ where $\delta\in [0, \delta_0]$.
\end{proof}

\bigskip

\noindent\textbf{Case~2.} Let $\theta_1 \geq 0, \theta_2 >0$ and $\beta \geq 0$. Consider the following FBDSDEJ:
\begin{eqnarray}\label{eq:3.26}
\left\{
\begin{array}{ll}
dy_{t}=[  \alpha \,  b\left(  t,\upsilon_{t}\right)  -\left(
1-\alpha\right)  \, \theta_{2} \,  Y_{t}+\widetilde{b}_{0}(t)
] \, dt-z_{t}\, \overleftarrow{dB}_{t}
 \\  \hspace{2in} + \, [  \alpha \,  \sigma\left(  t,\upsilon_{t}\right) -\left(
1-\alpha \, \right)  \theta_{2} \,  Z_{t}+\widetilde{\sigma}_{0}(t)
] \,  dW_{t}
\\  \hspace{0.9in} +\int_{\Theta} \, [  \alpha \, \varphi\left(  t,\upsilon_{t},\rho\right)
-\left(  1-\alpha\right)  \theta_{2} \,  k_{t}\left(  \rho\right)
+\varphi_{0}\left(  t,\rho\right)  ]\,  \widetilde{N}\left(
d\rho,dt\right), \\ \\
dY_{t}=[  \alpha \,  f\left(  t,\upsilon_{t}\right)  +\widetilde{f}_{0}(t)  ]  \, dt+[  \alpha \,  g\left(
t,\upsilon_{t}\right)  +\widetilde{g}_{0}(t)  ]
\, \overleftarrow{dB}_{t} \\  \hspace{2.8in} +\,  Z_{t}\, dW_{t}
+\int_{\Theta}k_{t}\left(  \rho\right)  \widetilde{N}(d\rho,dt), \\ \\
y_{0}=x,Y_{T}=\alpha \,  h\left(  y_{T}\right)  +\phi.
\end{array}
\right.
\end{eqnarray}

When $\alpha=0$, this FBDSDEJ is uniquely solvable, as shown in Case~1, but now by handling its second equation first.
On the other hand, when $\alpha=1$, we realize that the
existence of the solution of (\ref{eq:3.26}) implies that of (\ref{eq:3.2}).
\begin{lemma}\label{lem:final-lemma}
Assume that $\theta_1 \geq 0, \theta_2 >0$ and $\beta \geq 0$. Assume also {\rm(A1)-(A4)}. Then there exists
a positive constant $\delta_{0}$\ such that if, a priori, for each element $\phi$ of
$L^{2}\left(  \Omega,\mathcal{F}_{T},\mathbb{P};H\right)  $\ and $\left(
\widetilde{b}_{0},\widetilde{f}_{0},\widetilde{\sigma}_{0},\widetilde{g}_{0},\varphi
_{0}\right)  \in\mathfrak{M}^{2}$, FBDSDEJ~(\ref{eq:3.26}) is uniquely solvable for some
$\alpha_{0}\in\left[  0,1\right)  $, then for each elements $\alpha \in \left[
\alpha_{0},\alpha_{0}+\delta_{0}\right]  $, $\phi\in L^{2}\left(
\Omega,\mathcal{F}_{T},\mathbb{P};H\right) $ \ and $\left(  \widetilde{b}%
_{0},\widetilde{f}_{0},\widetilde{\sigma}_{0},\widetilde{g}_{0},\varphi_{0}\right)
\in\mathfrak{M}^{2}$, this system (\ref{eq:3.26}) is also uniquely solvable in $\mathfrak{M}^{2}$.
\end{lemma}
\begin{proof}  Assume that, for each $\phi\in L^{2}\left(  \Omega,\mathcal{F}%
_{T},\mathbb{P};\mathbb{R}^{m}\right)  $\ and $(  \widetilde{b}_{0},\widetilde{f}
_{0}, \widetilde{\sigma}_{0},\widetilde{g}_{0},\varphi_{0})  \in\mathfrak{M}^{2}$,
there exists a unique solution of FBDSDEJ~(\ref{eq:3.3}) for $\alpha=\alpha_{0}$. Then, for each element
$\overline{\upsilon}_{t}:=\left(  \overline{y}_{t},\overline{Y}_{t},\overline{z}_{t},\overline{Z}
_{t},\overline{k}_{t}\right)$ of $\mathfrak{M}^{2}$, there exists a unique
$\upsilon:=\left(  y,Y,z,Z,k\right)$ in $\mathfrak{M}^{2}$\ satisfying the following FBDSDEJ:
\begin{eqnarray}\label{eq:3.28}
\left\{
\begin{array}{ll}%
dy_{t}=[  \alpha_{0}\, b\left(  t,\upsilon_{t}\right) -(1-\alpha_{0})\, \theta_{2}\, Y_{t} +\delta \left( b(  t,\overline{\upsilon}_{t})+\theta_{2}\,  \overline{Y}_{t}\right) +\widetilde{b}_{0}\left(
t\right)  ] \, dt \\ \hspace{0.75cm}
+\, [  \alpha_{0}\, \sigma\left(  t,\upsilon_{t}\right) -(1-\alpha_{0})\, \theta_{2}\, Z_{t} +\delta  \left(  \sigma(  t,\overline{\upsilon}_{t})+\theta_{2}\,  \overline{Z}_{t} \right)   +\widetilde{\sigma}_{0}\left(
t\right)  ] \, dW_{t}  \\ \hspace{0.75cm}
-\, z_{t}\,\overleftarrow{dB}_{t}  +\int_{\Theta}\, [  \alpha_{0}\, \varphi\left(  t,\upsilon_{t}
,\rho\right)  -(1-\alpha_{0})\, \theta_{2}\, k_{t}(\rho)\\ \hspace{3.2cm}
+\, \delta  \left(  \varphi(  t,\overline{\upsilon}_{t},\rho)+\theta_{2}\,  \overline{k}_{t}(\rho)\right)
+\varphi_{0}\left(  t,\rho\right)  ]  \, \widetilde{N}\left(  d\rho,dt\right),\\ \\
dY_{t}=\left[  \alpha_{0}\, f\left(  t,\upsilon_{t}\right) +\delta    f\left(
t,\overline{\upsilon}_{t}\right)   +\widetilde{f}
_{0}\left(  t\right)  \right]  dt + Z_{t}\, dW_{t} \\ \hspace{1cm} +\left[  \alpha_{0}\, g\left(  t,\upsilon_{t}\right) +\delta \,  g\left(
t,\overline{\upsilon}_{t}\right)   +\widetilde{g}
_{0}\left(  t\right)  \right] \, \overleftarrow{dB}_{t}
+\int_{\Theta}k_{t}\left(  \rho\right)  \widetilde{N}\left(
d\rho,dt\right) ,
\\ \\
y_{0}=x , \, Y_{T}=\alpha_{0}\, h \left(  y_{T}\right)  +\delta \,  h\left(  \overline{y}
_{T}\right)  +\phi .
\end{array}
\right.
\end{eqnarray}

We consider the mapping $I_{\alpha_{0}+\delta}$ as defined in the proof of Lemma~\ref{Lemma: 3.5}, and we continue with the same notations proceeding
system~(\ref{eq:new-number}). Then, we apply It\^{o}'s formula to compute $\langle \widehat{y}_s,\widehat{Y}_s\rangle
$\ over $\left[  0,T\right]$, take the expectation, and use hypothesis (A3) to get eventually
\begin{eqnarray*}\label{eq:3.29}
&& \hspace{-0.5cm}\mathbb{E} \left[  \left\langle \widehat{y}_{T},\alpha_{0}\, \widehat{h}\left(
y_{T}\right)+\delta \,   \widehat{h}\left(
\overline{y}_{T}\right)\right\rangle \right]=\, \alpha_{0}\, \mathbb{E} \left[  \int_{0}^{T}\left\langle \widehat{A}\left(
s,\upsilon_{s}\right)
,\widehat{\upsilon}_{s}\right\rangle ds\right]\nonumber \\
&&
 + \, \delta \, \mathbb{E} \left[  \int_{0}^{T}\left\langle \widehat{y}_{s}\, , \widehat{f}\left(  s,\overline{\upsilon}_{s}\right)  \right\rangle
ds\right]+ \delta \, \mathbb{E} \left[
\int_{0}^{T}\left\langle \widehat{z}_{s} \, , \widehat{g}\left(
s,\overline{\upsilon}_{s}\right) \right\rangle ds\right]\nonumber \\
&& +\, \mathbb{E} \left[  \int_{0}^{T}\left\langle \widehat{Y}_{s} \, , -\left(  1-\alpha_{0}\right)  \theta_{2}\, \widehat{Y}_{s} +\delta \,  \left( \widehat{b%
}\left(  s,\overline{\upsilon}_{s}\right)+\theta_{2}\, \widehat{\overline{Y}}_{s} \right) \right\rangle ds\right] \nonumber \\
&&  +\, \mathbb{E} \left[  \int_{0}^{T}\left\langle -\left(  1-\alpha_{0}\right)  \theta_{2}\, \widehat{Z}_{s} +\delta \,  \left( \widehat{\sigma%
}\left(  s,\overline{\upsilon}_{s}\right)+\theta_{2}\, \widehat{\overline{Z}}_{s} \right),\widehat{Z}_{s} \right\rangle ds\right] \nonumber \\
&&  +\, \mathbb{E} \left[  \int_{0}^{T}\int_{\Theta
}\left\langle -\left(  1-\alpha_{0}\right)  \theta_{2}\, \widehat{k}_{s}(\rho)+  \delta   \left( \widehat{\varphi%
}\left(  s,\overline{\upsilon}_{s},\rho\right)+\theta_{2}\, \widehat{\overline{k}}_{s}(\rho) \right),\widehat{k}_{s}(\rho) \right\rangle \Pi(d\rho)ds\right].
\end{eqnarray*}
Thus, by applying  (A1), (A2) and the Cauchy-Schwarz inequality (while recalling that 
$\alpha_0 <1$), we obtain
\begin{eqnarray*}
&& \hspace{-0.6cm} \alpha_{0} \, \beta\, \mathbb{E} \left[  |\widehat{y}_{T}|^{2} \right] + (1-\alpha_{0})\, \theta_{2}\, \mathbb{E} \left[  \int_{0}^{T}\left(
\vert \widehat{Y}_{s}\vert ^{2}+\Vert \widehat{Z}_{s}\Vert^{2}+||| \widehat{k}_{s}||| ^{2}\right)  ds\right]    \\
&&  + \,\alpha_{0}\, \theta_{2}\, \mathbb{E} \left[  \int_{0}^{T}\left(
\vert \widehat{Y}_{s}\vert ^{2}+\Vert \widehat{Z}_{s}\Vert^{2}+||| \widehat{k}_{s}||| ^{2}\right)  ds\right]+ \alpha_{0} \, \theta
_{1}\, \mathbb{E} \left[  \int_{0}^{T}\left(  \left\vert \widehat{y}_{s}\right\vert ^{2}+\left\Vert \widehat{z}%
_{s}\right\Vert^{2}\right)  ds\right] \\
&& \hspace{-0.25cm} \leq \delta  \left(\mathbb{E} \left[  |\widehat{y}_{T}|^{2}\right] +c\, \mathbb{E} \left[  |\widehat{\overline{y}}_{T}|^{2}\right]\right) + \delta\, \mathbb{E} \left[  \int_{0}^{T} \left(  \left\vert \widehat{y}_{s}\right\vert ^{2} + \| \widehat{z}_{s}\| ^{2} \right) ds\right]  + \delta \, C\, \mathbb{E} \left[  \int_{0}^{T}\| \widehat{\overline{\upsilon}}_{s}\|^{2}ds \right] \\
&& \hspace{3cm}   + \, \delta \left( 1+\theta_2\right)  \mathbb{E} \left[  \int_{0}^{T} \left( \vert \widehat{Y}_{s}\vert ^{2} + \| \widehat{Z}_{s}\|^{2}
+ ||| \widehat{k}_{s}||| ^{2} \right)  ds \right]   \\
&& \hspace{2.2in} +\, \delta \, \theta_{2}\, \mathbb{E} \left[  \int_{0}^{T} \left(\vert \widehat{\overline{Y}}_{s}\vert ^{2} + \| \widehat{\overline{Z}}_{s}\| ^{2} +
||| \widehat{\overline{k}}_{s}||| ^{2} \right) ds \right] .
\end{eqnarray*}

We conclude that
\begin{eqnarray}\label{eq:3.30}
&& \hspace{-1.5cm}\theta_{2}\, \mathbb{E} \left[\int
_{0}^{T}\left( \vert \widehat{Y}_{s}\vert ^{2}+\Vert \widehat{Z}_{s}\Vert^{2}+||| \widehat{k}_{s}||| ^{2}\right)  ds \right]
\nonumber \\ && \hspace{-0.5cm}
\leq \delta  \, \mathbb{E} \left[  \left\vert \widehat{y}_{T}\right\vert ^{2}\right]  +\delta \, C\, \mathbb{E} \left[  \vert \widehat{\overline{y}
}_{T} \vert ^{2}\right]  +\, \delta \, C \,\mathbb{E} \left[\int_{0}^{T}\left( \Vert \widehat{\upsilon}%
_{s}\Vert^{2}+\Vert \widehat{\overline{\upsilon}}
_{s}\Vert^{2}\right)  ds\right],
\end{eqnarray}
after neglecting $\alpha_{0}\, \theta_{1}\, \mathbb{E} \, \big[  \int_{0}^{T}\left(  \left\vert
\widehat{y}_{s}\right\vert ^{2}+\left\Vert \widehat{z}_{s}\right\Vert
^{2}\right)  ds \big]$ and $\alpha_{0}\, \beta  \; \mathbb{E} \left[
\left\vert \widehat{y}_{T}\right\vert ^{2}\right]$ which contain  $\alpha_{0}$.
It turns out that we have to find estimates for $\mathbb{E}\, [ \left\vert \widehat{y}%
_{T}\right\vert ^{2}]$ and $\mathbb{E} \, \big[  \int_{0}^{T}\left(  \left\vert
\widehat{y}_{s}\right\vert ^{2}+\left\Vert \widehat{z}_{s}\right\Vert
^{2}\right)  ds \big]$. So, we apply It\^{o}'s formula to $\left\vert \widehat{y}_{s}\right\vert ^{2}$ over the interval
$\left[  0,t\right]$ to deduce that
\begin{eqnarray}\label{eq:3.31}
&& \hspace{-0.85cm}\mathbb{E} \left[  \left\vert \widehat{y}_{t}\right\vert ^{2}\right]
+\mathbb{E} \left[  \int_{0}^{t}\left\Vert \widehat{z}_{s}\right\Vert
^{2}ds\right]  \nonumber \\
&& \hspace{-0.5cm} \leq \, 2\, \mathbb{E} \left[  \int_{0}^{t}\left\vert \widehat{y}_{s}\right\vert
\cdot\left(  \alpha_{0}\, \left\vert \widehat{b}\left(  s,\upsilon
_{s}\right)  \right\vert +\left(  1-\alpha_{0}\right)  \theta_{2}\, \vert
\widehat{Y}_{s}\vert
+ \delta \left\vert \widehat{b}\left(
s,\overline{\upsilon}_{s}\right)  \right\vert +\delta \, \theta_{2}\, \vert
\widehat{\overline{Y}}_{s}\vert \right)  ds\right]
\nonumber \\
&& \hspace{2cm} +\, \left(\frac{1+\gamma}{2\gamma}\right) \, \alpha_{0}^{2}\, \mathbb{E} \left[  \int_{0}^{t}\left\vert \widehat{\sigma}\left(  s,\upsilon_{s}\right)  \right\vert
^{2}ds\right]
\nonumber\\
&& +\, 3\left(  \frac{1+\gamma}{1-\gamma
}\right) \mathbb{E} \left[  \int_{0}^{t} \left(  \left(  1-\alpha_{0}\right)^{2}\theta_{2}^{2} \, \Vert \widehat
{Z}_{s}\Vert^{2}+\delta^{2} \left\vert \widehat{\sigma}\left(
s,\overline{\upsilon}_{s}\right)  \right\vert ^{2}
+ \delta^{2}\, \theta_{2}^{2}\, \Vert \widehat{\overline{Z}}_{s}\Vert^{2}\right)  ds\right]
\nonumber \\
&& \hspace{2cm} +\, \left(\frac{1+\gamma}{2\gamma}\right) \alpha_{0}^{2} \, \mathbb{E} \left[  \int_{0}^{t} \int_{\Theta}\left\vert \widehat{\varphi}\left(  s,\upsilon_{s},\rho\right)  \right\vert
^{2}\Pi_{s}(d\rho) \, ds\right]
\nonumber\\
 && +\, 3\left(  \frac{1+\gamma}{1-\gamma}\right) \, \mathbb{E} \left[  \int_{0}^{t} \int_{\Theta}\left(  \left(  1-\alpha_{0}\right)^{2}\theta_{2}^{2} \,| \widehat{k}_{s}(\rho)|^{2}+\delta^{2}\, \left\vert \widehat{\varphi}\left(
s,\overline{\upsilon}_{s},\rho\right)  \right\vert ^{2}\right. \right. \nonumber \\
&& \hspace{3.2in} + \left.  \left. \delta^{2}\, \theta_{2}^{2}\, | \widehat{\overline{k}}_{s}(\rho)|^{2}\right)\Pi_{s}(d\rho) \, ds \right]
\nonumber \\
&& =:I'_{1}(t)+I'_{2}(t)+I'_{3}(t)+I'_{4}(t)+I'_{5}(t).
\end{eqnarray}

As we have previously shown in the proof of Lemma~\ref{Lemma: 3.5}, it is possible to derive estimates for $I'_{1}(t)-I'_{5}(t)$. By applying these estimates in (\ref{eq:3.31}) and utilizing Gronwall's inequality, we can obtain
\begin{eqnarray}\label{eq:3.38}
&& \hspace{-1cm} \mathbb{E} \left[  \left\vert \widehat{y}_{T}\right\vert^{2}\right] +\mathbb{E} \left[  \int_{0}^{T}  \left(  \left\vert
\widehat{y}_{s}\right\vert ^{2}+\left\Vert \widehat{z}_{s}\right\Vert
^{2}\right)
ds\right]  \nonumber \\ && \hspace{1cm}
\leq  C\, \mathbb{E} \left[  \int_{0}^{T}\left(
\vert \widehat{Y}_{s}\vert ^{2}+\Vert \widehat{Z}_{s}\Vert
^{2}+||| \widehat{k}_{s}|||
^{2}+\delta \, \Vert \widehat{\overline{\upsilon}}_{s}\Vert ^{2}\right)  ds\right] .
\end{eqnarray}

Now, with these crucial inequalities (\ref{eq:3.30}) and (\ref{eq:3.38}), we can proceed with the same reasoning as in the proof of Lemma~\ref{Lemma: 3.5} to deduce that
\begin{eqnarray*}
&& \hspace{-0.5cm} \mathbb{E} \left[  \left\vert \widehat{y}%
_{T}\right\vert ^{2}\right]  + \mathbb{E} \left[  \int_{0}^{T}\left\Vert \widehat{\upsilon}_{s}\right\Vert^{2}ds\right]  \leq\frac{1}{2}\left(  \mathbb{E} \left[  \vert \widehat{\overline{y}}_{T}\vert ^{2}\right]
+ \mathbb{E} \left[  \int_{0}^{T}\Vert \widehat
{\overline{\upsilon}}_{s}\Vert^{2}\, ds\right] \right) .
\end{eqnarray*}
Hence,  $I_{\alpha_{0}+\delta}$ is a contraction on  $\mathfrak{M}^{2}\times L^{2}(\Omega,\mathcal{F}_{T},\mathbb{P};\mathbb{R}^{n}),$ for all fixed $\delta \in[0,\delta_{0}]$ and some small $\delta_0$.  Consequently, $I_{\alpha_{0}+\delta}$ has a unique fixed point $\upsilon=\left(
y,Y,z,Z,k\right)\in \mathfrak{M}^{2},$\ which is evidently the unique solution of (\ref{eq:3.26}) for
$\alpha=\alpha_{0}+\delta, \delta \in [0, \delta_0]$.
\end{proof}

\bigskip

We are now prepared to finalize the proof of Theorem~\ref{Propo 3.4}.

\medskip

\noindent\begin{proof}[Proof completion of Theorem~\ref{Propo 3.4}]
In Case~1, we have already established in its proof that for every $\phi \in L^{2}\left(
\Omega,\mathcal{F}_{T},\mathbb{P};H \right)$ and $\left(  \widetilde{b}
_{0},\widetilde{f}_{0},\widetilde{\sigma}_{0},\widetilde{g}_{0},\varphi_{0}\right)
\in\mathfrak{M}^{2}$, the FBDSDEJ~(\ref{eq:3.3}) has a unique solution as $\alpha=0.$ Thus, by applying Lemma~\ref{Lemma: 3.5}, we can find $\delta_{0} > 0$ that  depends solely on the constants
$c,\gamma, \theta_{1}, \theta_{2}$, and $T$\ such that for any $\delta
\in\left[  0,\delta_{0}\right], \phi\in L^{2}\left(
\Omega,\mathcal{F}_{T},\mathbb{P}; H \right)$,  and $\left(  \widetilde{b}
_{0},\widetilde{f}_{0},\widetilde{\sigma}_{0},\widetilde{g}_{0},\varphi_{0}\right)
\in\mathfrak{M}^{2}$, the FBDSDEJ~(\ref{eq:3.3}) has a unique solution for $\alpha=\delta$. Now, since $\delta_{0}$ depends only on the aforementioned  constants, we can repeat this process $N$ times with ${1\leq N\delta_{0}<1+\delta_{0}}$. Therefore, for $\alpha=1$ specifically with $\left(  \widetilde{b}_{0},\widetilde{f}_{0}
,\widetilde{\sigma}_{0},\widetilde{g}_{0},\varphi_{0}\right)  \equiv0$ and $\phi\equiv
0$, we conclude that the FBDSDEJ~(\ref{eq:3.1}) attains a unique solution in 
$\mathfrak{M}^{2}$.

\medskip

In Case 2, we observe that for each $\phi \in L^{2}\left(
\Omega,\mathcal{F}_{T},\mathbb{P}; H \right)$, $\left(  \tilde{b}
_{0},\widetilde{f}_{0},\widetilde{\sigma}_{0},\widetilde{g}_{0},\varphi_{0}\right)\in\mathfrak{M}^{2}$, the FBDSDEJ~(\ref{eq:3.26}) has a unique solution as $\alpha=0$. Consequently,
Lemma~\ref{lem:final-lemma} implies that there exists a positive constant $\delta_{0}$ that depends solely on the constants
$c,\gamma, \theta_{1}, \theta_{2}$, and $T$\ such that, for any $\delta \in \left[  0,\delta_{0}\right], \phi\in L^{2}\left(
\Omega,\mathcal{F}_{T},\mathbb{P}; H \right)$, and $\left(  \widetilde{b}
_{0},\widetilde{f}_{0},\widetilde{\sigma}_{0},\widetilde{g}_{0},\varphi_{0}\right)
\in\mathfrak{M}^{2}$, the FBDSDEJ~(\ref{eq:3.26}) has a unique solution for $\alpha=\delta$. Hence, similar to the previous case, we can deduce that FBDSDEJ~(\ref{eq:3.1}) has a unique solution in $\mathfrak{M}^{2}$.
\end{proof}

\begin{example}\label{ex:ex1}
Let $E$ and $H$ be two real separable Hilbert spaces. Assume that $B$ and $W$ are cylindrical Wiener process on $E$. We consider the following system on $H$:
\begin{equation}\label{eq:3.9}
\left\{
\begin{array}{ll}
dy_{t}=-Y_{t} \, dt+\frac{1}{4} \left( z_{t} - Z_{t}\right) dW_{t}-z_{t}\, \overleftarrow{d B}_{t} - \frac{1}{4}  \, k_{t}\left(\rho\right)
\widetilde{N}(d\rho,dt), &  \\
dY_{t}= - y_{t} \, dt
 -  \left( z_{t} + \frac{1}{4}\,  Z_t \right) \overleftarrow{d B}_{t}+ Z_{t}\, dW_{t} + \int_{\Theta} k_{t}\left(  \rho\right)  \widetilde
{N}(d\rho,dt)  , &  \\
y_{0}=x \; (\in H),\quad Y_{T}=y_T .
&
\end{array}%
\right.
\end{equation}

In order to establish a connection between this system and FBDSDEJ~(\ref{eq:3.1}), we define, for $t\in [0,T]$, the following mappings:
\begin{equation*}
\begin{array}{ll}
b\left( t,y_{t},Y_{t},z_{t},Z_{t},k_{t}\right)= - Y_{t}, &  \\
\sigma\left( t,y_{t},Y_{t},z_{t},Z_{t},k_{t}\right)=\frac{1}{4} \left( z_{t} - Z_{t}\right), &  \\
\varphi\left( t,y_{t},Y_{t},z_{t},Z_{t},k_{t}, \rho\right)=- \frac{1}{4}  \, k_{t}\left(\rho\right), &  \\
f\left( t,y_{t},Y_{t},z_{t},Z_{t},k_{t}\right)= y_{t}, &  \\
g\left( t,y_{t},Y_{t},z_{t},Z_{t},k_{t}\right)=z_{t} + \frac{1}{4}\,  Z_t , &  \\
h\left( y_{T}\right) = y_{T}.
\end{array}%
\end{equation*}

It is evident that the assumptions {\rm(A1)--(A4)} are satisfied for this example with the constants: $c=1, \gamma=\frac{1}{4} \, , \,  \theta_{1}=\theta _{2}=\frac{1}{4} \,, $ and $\beta=1$. Therefore,
Theorem~\ref{Propo 3.4-1} and Theorem~\ref{Propo 3.4} guarantee the existence of the unique solution of FBDSDEJ~(\ref{eq:3.9}).
\end{example}

\begin{example}\label{ex:ex2}
Let us consider the following FBDSDEJ on the space $H=\mathbb{R}$:
\begin{equation}\label{eq:3.10-ex}
\left\{
\begin{array}{ll}
dy_{t}= Y_{t} \, dt -z_{t} \, \overleftarrow{d B}_{t} , &  \\
dY_{t}=- \, y_{t} \, dt-z_{t}\, \overleftarrow{d B}_{t}+Z_{t}\, dW_{t}  + \int_{[0,1]}k_{t}\left(  \rho\right)  \widetilde{N}(d\rho,dt) , &  \\
y_{0}=0,\hspace{0.5cm}Y_{T}=-\, y_{T}  , &
\end{array}%
\right.
\end{equation}
where $T= \frac{3 \pi}{4}$, and $B$ and $W$ are $1$-dimensional Brownian motions.

Using the notation of assumption {\rm (A1)}, for $v=\left( y,Y,z,Z,k\right),$  we have
\begin{equation*}
\begin{array}{ll}
A\left( t,v_t\right) =\left( - \,
y_{t}  ,
Y_{t} ,-z_{t},0, 0\right)  . &
\end{array}%
\end{equation*}
Since
\begin{equation*}
\left\langle A\left( t,v_t^{1}\right) -A\left(
t,v_t^{2}\right) ,v_t^{1}-v_t^{2} \, \right\rangle  =- \left(
\Delta y_{t}\right) ^{2} + \left(  \Delta Y_{t} \right) ^{2}-%
 \left( \Delta z_{t}\right) ^{2},
\end{equation*}
we observe that {\rm (A1)} is not fulfilled.
As a result, (\ref{eq:3.10-ex}) may not possess a unique solution. In fact, $(\sin t, \cos t, 0,0,0)$ is  a solution of (\ref{eq:3.10-ex}) in addition to the trivial solution $(y_t,Y_t,z_t,Z_t,k_t)=(0,0,0,0,0)$.
\end{example}

\fussy

\end{document}